# FFT, FMM, OR MULTIGRID? A COMPARATIVE STUDY OF STATE-OF-THE-ART POISSON SOLVERS FOR UNIFORM AND NON-UNIFORM GRIDS IN THE UNIT CUBE

AMIR GHOLAMI*, DHAIRYA MALHOTRA*, HARI SUNDAR*, AND GEORGE BIROS*

**Abstract.** From molecular dynamics and quantum chemistry, to plasma physics and computational astrophysics, Poisson solvers in the unit cube are used in many applications in computational science and engineering. In this work, we benchmark and discuss the performance of the scalable methods for the Poisson problem which are used widely in practice: the Fast Fourier Transform (FFT), the Fast Multipole Method (FMM), the geometric multigrid (GMG) and algebraic multigrid (AMG). Our focus is on solvers support high-order, highly non-uniform discretizations, but for reference we compare with solvers specialized for problems on regular grids. So, we include FFT, since it is a very popular algorithm for several practical applications, and the finite element variant of HPGMG, a high-performance geometric multigrid benchmark. In total we compare five different codes, three of which are developed in our group. Our FFT, GMG and FMM are parallel solvers that use high-order approximation schemes for Poisson problems with continuous forcing functions (the source or right-hand side). Our FFT code is based on the FFTW for single node parallelism. The AMG code is from the Trilinos library from the Sandia National Laboratory. Our geometric multigrid and our FMM support octree based mesh refinement, variable coefficients, and enable highly non-uniform discretizations. The GMG, actually also supports complex (non-cubic) geometries using a forest of octrees.

We examine and report results for weak scaling, strong scaling, and time to solution for uniform and highly refined grids. We present results on the Stampede system at the Texas Advanced Computing Center and on the Titan system at the Oak Ridge National Laboratory. In our largest test case, we solved a problem with 600 billion unknowns on 229,379 cores of Titan. Overall, all methods scale quite well to these problem sizes. We have tested all of the methods with different source functions (the right hand side in the Poisson problem). Our results indicate that FFT is the method of choice for smooth source functions that require uniform resolution. However, FFT loses its performance advantage when the source function has highly localized features like internal sharp layers. FMM and GMG considerably outperform FFT for those cases. The distinction between FMM and GMG is less pronounced and is sensitive to the quality (from a performance point of view) of the underlying implementations. In most cases, high-order accurate versions of GMG and FMM significantly outperform their low-order accurate counterparts.

**1. Introduction.** The need for large scale parallel solvers for elliptic partial differential equations (PDEs) pervades across a spectrum of problems with resolution requirements that cannot be accommodated on current systems. Several research groups are working on technologies that scale to trillions of unknowns and billions of cores. To illustrate some of the issues in scaling such solvers and to provide a (non-exhaustive) snapshot of current technologies, we conduct an experimental study of solving a simple model elliptic PDE and compare several state of the art methodologies.

We restrict our attention to the following model problem: given $f$, a smooth and periodic function in the unit cube, we wish to find $u$ (also smooth and periodic in the unit cube) such that

$$-\Delta u = f, \qquad (1.1)$$

*Institute for Computational Engineering and Sciences, The University of Texas at Austin, Austin, TX 78712. (i.amirgh@gmail.com, dhairya.malhotra@gmail.com, hsundar@gmail.com, gbiros@acm.org).

where $\Delta$ is the Laplace operator. This is also known as the constant-coefficient Poisson problem. It encapsulates many of the difficulties in solving elliptic partial differential equations (PDEs). We chose this problem because all four methods can address it in an algorithmically optimal way. Algorithms for solving this problem, also known as "*Poisson solvers*" find applications in astrophysics, chemistry, mechanics, electromagnetics, statistics, and image processing, to name a few. Vendors like Intel and NVIDIA provide Poisson solvers in their math libraries. Examples of scientific computing libraries that provide Poisson solvers include PETSc [3], Trilinos [24], deal.II [4], and MATLAB.

Poisson solvers must scale to trillions of unknowns. Example of methods that scale well are the FFT (based on spectral discretizations)[1], the Fast Multipole Method (based on discretizing the integral equation reformulation of (1.1), and multigrid methods (for stencil-based discretizations). Other scalable methods include domain decomposition and wavelet transforms, which will not be discussed here. Despite the existence of many different Poisson solvers there has been little work in *directly* benchmarking the computational efficiency of these methods, in particular for the case of non-uniform discretizations. Such benchmarking is quite typical in other scientific computing areas (e.g., sorting, matrix computations, and graph partitioning).

*Methodology and contributions.* In this paper we benchmark five state-of-the-art algorithms and implementations, three from our group and two from different groups:

The first solver is parallel FFT using the AccFFT which has been recently developed in our group [22]. AccFFT, built on top of FFTW, uses MPI and OpenMP, as well as novel communication schemes that makes it faster than similar libraries. We report comparisons with PFFT [36] and P3DFFT [33].

The second solver is the ML algebraic multigrid solver which is part of the Trilinos library developed and maintained by the U.S. Department of Energy [21, 24]. In our runs, we use ML with MPI. It is one of the most scalable, general purpose codes available. It can handle much more complex problems than the one we consider here.

The third solver is an in-house Geometric Multigrid scheme that uses continuous Galerkin discretizations on octree meshes. The low-order version of the code appeared in [39], but the scalable high-order results we report here are new, as well as the outline of the algorithm. The sequential algorithm is described in [41]. The library uses MPI and more details about the new algorithm will be presented in §2.

The fourth code is PvFMM, also an in-house novel parallel volume FMM that supports continuous as well as particle sources ($f$). It uses octree discretization using Chebyshev polynomials at each leaf node to represent $f$. PvFMM uses MPI and OpenMP [30].

The fifth code is HPGMG [1] and it is a high-performance computing benchmark code for regular grids with finite-element and finite-volume implementations. The finite-element implementation is more general as it supports variable coefficients and coordinate transformation Jacobians that resemble overset grids and non-uniform grids. For this reason we include only the finite-element version in our tests.

We compare these five methods on two different architectures, Stampede and Titan, and we discuss two main questions:
- Which method is faster? What matters the most is the wall-clock time to solution. That is, given $f$, we would like to evaluate $u$ (typically at a given number of points) to a specified accuracy. We consider two main cases, highly oscillatory fields, for which a regular grid is necessary, and highly localized fields for which adaptively refined meshes are expected to be more effective.
- How does the cost per unknown compare for the different methods given a fixed algebraic accuracy? This test focuses on constant in the complexity estimates which are functions of the problem size, the number of processors, the approximation order, and of course the implementation. We perform weak and strong scaling studies to directly compare the complexity estimates on specific architectures using the same problem size. We have scaled our runs up to 229,376 cores on Titan for problems with up to 600 billion unknowns. Our goal here is not to fit a detailed performance model to the runs, but rather to identify whether there are order-of-magnitude differences in the performance between these methods.

---

[1] FFT can be used also to diagonalize and invert stencil discretizations on uniform grids. We are not discussing this case here.

To our knowledge, such a benchmark at these problem sizes and number of cores is the first of its kind. We view it as a companion to existing theoretical complexity analysis. In addition to work complexity given a problem size, we consider the issue of work complexity given a target accuracy using both uniform and non-uniform grids—the "right" problem size is *not known* a priori. One reason such a study has not taken place is that the underlying technologies have not been available at this scale. Indeed, we are not aware of any other distributed-memory FMM codes that allow $f$ to be an arbitrary function (most existing codes only support sums of delta functions, also known as point-FMMs methods). Also, the only other scalable, high-order, multigrid scheme we know is that of Paul Fischer's group [29]. Both of our GMG and FMM codes support arbitrary order discretizations. In summary, we test weak and strong scalability of all of these methods and report time to solution, and setup time for different test cases.

Qualitatively, the results of our study for solving (1.1) can be summarized as follows: FFT is the method of choice for uniform discretizations even at large-core count problems. FMM and MG are the methods of choice in the presence of strongly-localized features. AMG scales well, but it is significantly slower, especially when including setup costs. Uniform grid, second-order discretizations (e.g. the 7-point Laplacian) end up being 1000× or more, slower than high-order schemes and the FFT for high-accuracy solutions. Even for low-accuracy solutions in non-uniform grids, second-order methods suffer. Third or higher order can offer significant speed-ups. Of course, these conclusions are valid only when the solution is smooth.

*Related work.* There is rich literature discussing the accuracy and scalability of FFT, FMM and multigrid but to our knowledge, little work has been done on directly comparing the efficiency of these schemes. In [2], the authors compare FMM, FFT, and GMG for particle summation (point sources) with periodic conditions. This is different from what we look here, which is a Poisson solver with continuous right-hand side, not point sources. For the particle mesh variant of the codes tested, there is a continuous source but it is uniform. Furthermore they only consider nearly uniform distributions of charges. This is fine for molecular dynamics, but it is of rather narrow scope. Here we consider highly non-uniform sources. The results FFT and GMG solvers are not considered separately but only as part of particle-mesh solvers that include many additional components that complicate the interpretation of the results for other applications. In [20], theoretical complexity estimates for the scaling of FFT and multigrid (for uniform grids) are provided and their implications towards the design of exascale architectures is discussed. In [12], a similar study is carried for the FFT, along with experimental results on both CPU and hybrid systems. An interesting performance model is introduced that accounts for both intra-node and inter-node communication costs. In [5], the authors consider complexity estimates that account for low-level hardware details and consider the viability of different applications including FFT, and matrix vector multiplications, and molecular dynamics simulations. In [19], the authors discuss the scalability of algebraic multigrid (on uniform grids), provide performance models, and conduct an experimental scalability study on up to 65,536 cores. A perspective on scalability is given in [46]. In [49], the authors discuss the scalability of a point FMM code to exascale architectures and provide scalability results up to 32,768 cores and 40 billion unknowns. In our group we have worked on scalable geometric multigrid methods and their comparison to algebraic multigrid schemes [39] (but only for low-order discretizations), as well as parallel FMM schemes based on the kernel-independent variant of the FMM kernels [26, 48]. All of these studies are critical in understanding the scalability of the schemes and helping co-design the next architectures. We consider our study as a companion to these works as it provides experimental data that can be further analyzed using performance models. Also, except for the work in our group, others have only considered uniformly refined grids (and low-order discretizations for the multigrid). Here we consider all cases: uniform and refined grids, low and high-order discretizations, and four major algorithms.

Other scalable approaches to solving the Poisson problem include hybrid domain decomposition methods [29]. A very efficient Poisson solver is based on a non-iterative domain decomposition method [31] using a low-order approximation scheme. In [25], that solver was compared with a high order volume FMM. The FMM solver required $4\times$–$100\times$ fewer unknowns. Other works based on FFTs, tree codes and multigrid that are highly scalable (albeit for low-order, or point FMM only) include [14,32,34,35,38].

*Limitations.* Our study is limited to problems with constant coefficients on the unit cube, with periodic boundary conditions, and smooth solutions. AMG is the only method that is directly applicable to general geometries and problems with complex coefficients. GMG methods also can handle certain types of complex geometries. FMM can be used for constant-coefficient problems on arbitrary geometries but can also be extended to variable coefficients using volume integral equations. Complex geometries are also possible but the technology for such problems is not as developed as for algebraic multigrid. In our tests we only use periodic boundary conditions but the results for the FMM and GMG apply to Neumann and Dirichlet problems on the unit cube. For our high-order geometric multigrid we use a smoother that heuristically works well but we do not have supporting theory that shows that we have made the best possible choice (see §2.3). The solve phase for these codes has been optimized significantly but it doesn't mean that it can be further improved. In §4, we compare our GMG with the stencil-based HPGMG library [1] to give a measure of its performance relatively to a highly optimized (but less general code). Another limitation with respect performance is the use of heterogeneous architectures.

Comparing highly specialized codes with much general purpose codes (like the AMG solvers we use here) is problematic but, we think, informative as it provides some data for future developments of software that must be as efficient as possible if it is to be used in the million and billion-core systems. As we mentioned other methods like domain decomposition or hybrids like particle-in-cell methods are not discussed. One salient disadvantage of high-order methods is that when used with explicit time stepping schemes (for example, due to coupling to a transport equation) they can lead to extremely small time steps. Finally, FFTs, FMM, and Multigrid can be also discussed in other contexts (e.g., Ewald sums, particle-in-cell methods, signal analysis), which we do not do here.

*Outline of the paper.* In section §2, we summarize complexity estimates for FFT, Multigrid, and FMM. For FMM and geometric multigrid we provide some more detail since some components of the underlying algorithms are new. In §3, we summarize the experimental setup, platforms, and the choice of the right hand sides. In §4, we present and discuss the results of our experiments.

**Notation:** We use $p$ to denote the number of cores $q$ to denote the order of polynomial approximation for $f$ and $u$, $m$ the FMM approximation order for the far field, and $N$ the total number of unknowns.

**2. Methods.** Here we describe the basic algorithmic components of each method, and their overall complexity (setup and solve) and the solve complexity ($T_{\text{FMM}}, T_{\text{FFT}}, T_{\text{GMG}}$). In all of our results we assume that $N/p \gg 1$ and the complexity estimates are stated for an uncongested hypercube topology and uniform grids.

**2.1. The Fast Fourier Transform.** The Fast Fourier Transform is an algorithm for computing the Discrete Fourier Transform (DFT) of a signal in $\mathcal{O}(N \log N)$. Several efficient implementations of single-node FFT exist (for example, FFTW [18], Intel MKL's FFT [47], IBM's ESSL library [16], etc). In addition, several FFT algorithms have been proposed for distributed machines (for example see [17,27,36,42,43]). In computing 3D FFTs, the key challenge is in dividing the data across the processes. One option is slab decomposition in which the data is partitioned into $p$ slabs or slices, each containing $N/p$ samples. Such an approach however limits the maximum number of process to the number of slices in the input data ($p \leq \mathcal{O}(N^{1/3})$).

This bound can be increased to $\mathcal{O}(N^{2/3})$ by using a 2D decomposition, which is also referred to as *pencil decomposition*. Consider a 3D function of size $N = N_x \times N_y \times N_z$ whose FFT we wish to compute on $p$ processes. AccFFT maps the processes into a 2D grid of size $p_x \times p_y = p$, each having a block (or "pencil") of size $N_x/p_x \times N_y/p_y \times N_z$ (or an equivalent partition). In the limit of $p = N_x \times N_y = N^{2/3}$ processes, each process will get a one dimensional pencil of the data (of length $N_z$). With the pencil decomposition, each process has all the data in one direction (e.g. $z$ direction), and partial data of the other two directions ($x$ and $y$). Consequently, the 1D FFT along the $z$ direction can be performed independent of other processes.

Following the local 1D FFT (along the $z$-direction), each process exchanges its data using all-to-all communication with other processes in the same row ($p_x$) to form the transpose. After a process receives the data from other processes, it computes the 1D FFT along the next direction. This process is once again repeated in the column direction ($p_y$). The data layout in each stage is shown schematically in Figure 1.

The communication cost of distributed FFT is given by $\mathcal{O}(\frac{N}{\sigma(p)})$, where $\sigma(p)$ is the bisection bandwidth of the network; for a hypercube it is $p/2$ [33]. The total execution time of one 3D FFT can be approximated by

$$T_{\text{FFT}} = \mathcal{O}\left(\frac{N \log N}{p}\right) + \mathcal{O}\left(\frac{N}{p}\right).$$

Because we are comparing different solvers, it is important that each of the implementations are optimal. To address this concern, we compare AccFFT's performance with two other libraries in Table 1. The first library is P3DFFT, written by Dmitry

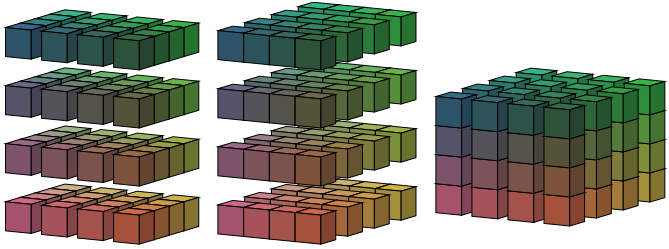

**Fig. 1:** *Data layout in different stages of computing a forward/inverse FFT using pencil decomposition for $p_x = p_y = 4$.*

|         | 512 Cores |       | 2048 Cores |       | 8192 Cores |       |
|---------|-----------|-------|------------|-------|------------|-------|
|         | *FFT*     | *Setup* | *FFT*    | *Setup* | *FFT*    | *Setup* |
| AccFFT  | 0.256     | 10.8  | 0.070      | 3.6   | 0.024      | 1.5   |
| P3DFFT  | 0.468     | 10.0  | 0.118      | 12.9  | 0.040      | 4.1   |
| PFFT    | 0.848     | 19.4  | 0.257      | 6.9   | 0.073      | 2.8   |

**Table 1:** Comparison of total time and setup time of AccFFT, P3DFFT, and PFFT libraries. We report timings for a single threaded, in-place real-to-complex FFT of size $1024^3$ on 512, 2048, and 8192 cores of TACC's Stampede platform. All the libraries were tested using `MEASURE` flag as well as other library specific options to achieve the best performance.

Pekurovsky [33]. P3DFFT is a robust parallel FFT library that has been successfully used in different applications and shown to have excelelent scalability. The second library is PFFT which has been recently released by Michael Pippig [36]. PFFT supports high dimensional FFTs as well as optimized pruned FFTs (with respect to time). It is designed to allow the user to compute transforms of more general data layouts. It has been tested up to 200K cores of BlueGene/Q, and shown to be scalable. Our AccFFT supports slab and pencil decompositions for CPU and GPU architectures, for real-to-complex, complex-to-complex, and complex-to-real transforms. It uses a series of novel schemes and in our tests results in an almost 2× speedup over P3DFFT and 3× speedup over PFFT as shown in Table 1. Details of all of the performance optimizations of AccFFT are beyond the scope of this paper. We refer to [22] for details on the method. We would like to remark that P3DFFT and PFFT have been designed to allow for more general features that are currently not available in AccFFT. In our tests we have used both P3DFFT and AccFFT.[2]

To solve the Poisson problem we compute the FFT of $f$ in (1.1), scale it by the corresponding inverse diagonal form of the Laplace operator (using a Hadamard product), and compute the inverse Fourier transform of the result. The overall complexity of solving the Poisson problem with FFT is the same order as $T_{\text{FFT}}$.

**2.2. The Fast Multipole Method.** The FMM was originally developed to speedup the solution of particle N-body problems by reducing the complexity from $\mathcal{O}(N^2)$ to $\mathcal{O}(N)$ [10, 11]. Solving Poisson's equation by computing volume potential is similar to an N-body problem, except that the summation over source points is replaced by an integral over the continuous source density. The work [15, 25] shows how one needs to modify the particle-FMM to obtain a volume FMM and this is the approach used in our implementation. Our contribution is that we have focused on the efficient implementation of this algorithm by considering efficient blocking and vectorization strategies for per core performance. We further optimized this method by adding support for the Intel Phi accelerator and extended it to a distributed memory method based on our hypercube communication algorithm [26]. For results in this paper, we do not use Phi since none of the other solvers can use it. Below we summarize the main components of our algorithm, which is described in full detail in [30].

---

[2] We have compiled all libraries with performance tuning on. This significantly increases the setup time. However, this cost can be amortized across multiple solves, so for most applications it is negligible.

We evaluate the solution to (1.1) as a convolution of the source density function $f$ with $K(x) = -(4\pi|x|)^{-1}$ (free space Green's function for Laplace equation). We accelerate computation of this convolution using FMM. The basic idea behind FMM is to construct a hierarchical decomposition of the computational domain using an octree. Then, the solution at each point $x$ can be evaluated by summing over contributions from all octants in the octree. This summation is split into near and far interactions:

$$u(x) = \sum_{B \in \text{Near}(x)} \int_B K(x-y)f(y) + \sum_{B \in \text{Far}(x)} \int_B K(x-y)f(y)$$

The near interactions (from $B \in \text{Near}(x)$) are computed through direct integration. The far interactions (from $B \in \text{Far}(x)$) are low-rank and can be approximated. In the following sections, we describe the tree construction and the different interactions in more detail.

**2.2.1. Octree Construction.** We partition the computational domain using an octree data structure. For each leaf octant $B$, we approximate the source density $f$ by Chebyshev polynomials of order $q$

$$f(y) = \sum_{\substack{i,j,k \geq 0 \\ i+j+k<q}} \alpha_{i,j,k} T_i(y_1) T_j(y_2) T_k(y_3)$$

where, $T_i$ is the Chebyshev polynomial of degree $i$ and $y$ is a point in the octant $B$. The absolute sum of the highest order coefficients is used as an estimate of the truncation error. This error estimate is used to refine adaptively until a specified error tolerance is achieved.

We also apply 2:1-balance constraint on our octree, that is we constrain the difference in depth of adjacent leaf octants to be at most one. To do this, we need to further subdivide some octants and is called 2:1-balance refinement [40].

**2.2.2. Far Interactions.** We define a source and a target octant to be well separated (or far) if they are at the same depth in the octree and are not adjacent. To compute far interactions we use two building blocks: multipole expansions and local expansions.[3] The multipole expansion approximates the potential of an octant far away from it. The local expansion approximates the potential within an octant due to sources far away from it. The interactions are then approximated by computing a multipole expansion (source-to-multipole, multipole-to-multipole) for the source octant, multipole-to-local (V-list) translation and then evaluating the local expansion (local-to-local, local-to-target) at the target octant. We use the kernel-independent variant of FMM in our implementation. The form of the multipole and local expansions and the V-list translation operator for this variant are discussed in detail in [48].

---
[3]We refer to multipole expansions and we use the term *"multipole order"* $m$ to denote the accuracy of the far field approximation, following the analytic FMM. However, in our implementation we use the kernel independent FMM and $m$ refers to the square root of the number of equivalent points. We used the terminology from the original FMM for the readers that are not familiar with the kernel independent FMM [48].

**2.2.3. Near Interactions.** The near interactions (source-to-target or U-list), between pairs of adjacent source and target octants, are computed using exact integration. Since these are singular and near-singular integrals (because the kernel has a singularity), it is not feasible to do this on the fly. We have to precompute integrals over Chebyshev basis functions, then the integral can be represented as a linear transformation:

$$u(x) = \int_B K(x-y)f(y) = \int_B K(x-y) \sum_{i,j,k} \alpha_{i,j,k} T_{i,j,k}(y) = \sum_{i,j,k} \alpha_{i,j,k} I_x$$

We precompute these integrals for each target point $x$, and for each interaction direction then construct interaction matrices. These interactions can then be evaluated through matrix-vector products. To limit the number of possible interactions directions and make this precomputation possible, we need the 2:1-balance constraint (discussed above).

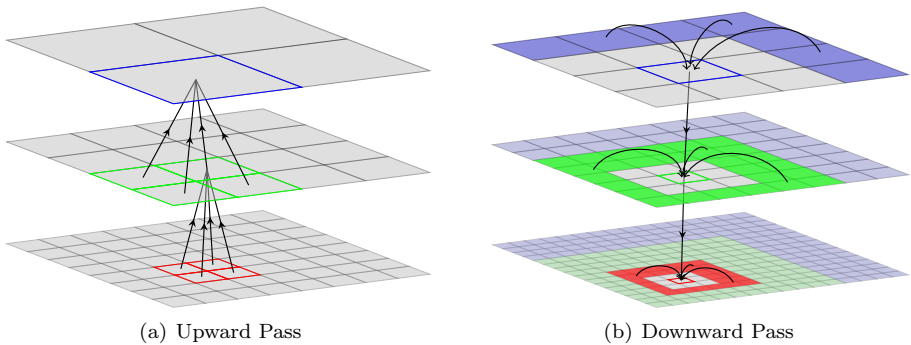

(a) Upward Pass    (b) Downward Pass

**Fig. 2:** *(a) Upward Pass: constructing multipole expansions. (b) Downward Pass: constructing local expansions, evaluating near interactions.*

**2.2.4. Summary of Volume FMM.** The overall algorithm for volume FMM can be summarized as follows:
- *Tree Construction*: Construct a piecewise Chebyshev approximation of the source density using octree based domain decomposition and perform 2:1-balance refinement.
- *Upward Pass:* For all leaf nodes apply source-to-multipole (S2M) translation to construct the multipole expansion from the Chebyshev approximation. For all non-leaf nodes apply multipole-to-multipole (M2M) translations in bottom-up order, to construct multipole expansion from that of its children, as shown in Fig 2(a).
- *Communication:* In the distributed memory implementation, we communicate the source density and the multipole expansions for ghost octants.
- *Downward Pass:* For all octree nodes, apply V-list and source-to-local (X-list) translations to construct the local expansion of each octant. In top-down order apply the local-to-local (L2L) translation to all nodes (Fig 2(b)). For all leaf octants, apply local-to-target (L2T), multipole-to-target (W-list) and source-to-target (U-list) translations to construct the final target potential as piecewise Chebyshev interpolation.

In our results, we report tree construction as the setup phase. The upward-pass, communication and downward-pass together constitute the evaluation phase.

**2.2.5. Parallel Fast Multipole Method.** Here we list the important features of the intra-node parallelism (accelerators, multithreading and vectorization) and distributed memory parallelism. A detailed discussion of these optimizations is beyond the scope of this paper. We refer the interested readers to [30] for a detailed discussion of these concepts.

*U,W,X-List Optimizations:* We group similar interactions, those with the same interaction matrix or related by a spatial symmetry relation into a single matrix-matrix product, evaluated efficiently through DGEMM.

*V-List Optimizations:* The V-list interactions involve computation of the Hadamard products, which has low computational intensity and is therefore bandwidth bound. We rearrange data and use spatial locality of V-list interactions to optimize cache utilization. This along with use of AVX and SSE vector intrinsics and OpenMP allowed us to achieve over 50% of peak performance for this operation on the Intel Sandy Bridge architecture.

*Distributed 2:1 Balance Refinement:* We developed a new distributed memory algorithm for 2:1 balance refinement, which is more robust for highly non-uniform octrees than our earlier implementation.

*Distributed Memory Parallelism:* We use Morton ordering to partition octants across processors during the tree construction. In the FMM evaluation, after the upward pass, we need to construct the local essential tree through a reduce-broadcast communication operation. For this, we use the hypercube communication scheme of [26].

**2.2.6. Complexity.** The cost of FMM evaluation is given by the number of interactions between the octree nodes weighted by the cost of each translation. The cost of each interaction depends on the multipole order $m$ and the order of Chebyshev polynomials $q$. Let $N_{oct}$ be the local octree nodes, $N_{leaf}$ the number of local leaf nodes. Also, let $N_U, N_V, N_W, N_X(= N_W)$ denote the number of interactions of each type U,V,W and X-list respectively. The overall cost is summarized in Table 2.

| Interaction Type | Computational Cost |
|---|---|
| S2M, L2T | $\mathcal{O}(N_{leaf} \times q^3 \times m^2)$ |
| M2M, L2L | $\mathcal{O}(N_{oct} \times m^2 \times m^2)$ |
| W-list, X-list | $\mathcal{O}(N_W \times m^2 \times q^3)$ |
| U-list | $\mathcal{O}(N_U \times q^3 \times q^3)$ |
| V-list | $\mathcal{O}(N_V \times m^3 + N_{oct} \times m^3 \log m)$ |

**Table 2:** *Computational cost for each interaction type.*

The communication cost for the hypercube communication scheme is discussed in detail in [26]. For an uncongested network, that work provides a worst case complexity which scales as $\mathcal{O}(N_s(q^3 + m^2)\sqrt{p})$, where $N_s$ is the maximum number of shared octants per processor. However, assuming that the messages are evenly distributed across processors in every stage of the hypercube communication process, we get a cost of $\mathcal{O}(N_s(q^3 + m^2) \log p)$. In our experiments with uniform octrees, the observed complexity agrees with this estimate. Also, since shared octants are near the boundary of the processor domains, we have $N_s \sim (N_{oct}/p)^{2/3}$, where $N_{oct}$ is total number of octants.

In the uniform octree case, there are $N_{oct} = N/q^3$ total octants, $N_U = 27N_{oct}$ and $N_V = 189N_{oct}$. Due to the large constant factors, the cost of U-list and V-list interactions dominate over other interactions and the overall cost is:

$$T_{\text{FMM}} = \mathcal{O}\left(q^3 \frac{N}{p}\right) + \mathcal{O}\left(\frac{m^3}{q^3} \frac{N}{p}\right) + \mathcal{O}\left(\left(\frac{N}{p}\right)^{2/3} q \log p\right).$$

In our results, we also report the time for setup (tree construction and 2:1 balance refinement). For tree construction, the cost of Chebyshev approximation for $\mathcal{O}(N/q^3)$ octants, distributed across $p$ processors is $\mathcal{O}(qN/p)$. During the adaptive refinement, we may also need to repeatedly redistribute octants across processors. This communication cost is data dependent and is difficult to analyze. For the 2:1 balance refinement, the cost is $\mathcal{O}\left(qN/p + (N \log N)/(pq^3)\right)$ assuming a hypercube interconnect.

**2.3. Geometric and Algebraic Multigrid.** Multigrid [6, 23, 45] is one of the most effective solvers for elliptic operators. It is algorithmically optimal and easy to implement and parallelize for uniform grids. Multigrid consists of two complimentary stages: *smoothing* and *coarse grid correction*. Smoothing involves the application of a (typically stationary) iterative solver to reduce (oscillatory) high-frequency errors. Coarse-grid correction involves transferring information to a coarser grid through *restriction*, solving a coarse-grid system of equations, and then transferring the solution back to the original grid through *prolongation* (interpolation). The coarse-grid correction eliminates (smooth) low-frequency errors. This approach can be applied recursively to obtain a multilevel system consisting of progressively coarser meshes. The multigrid method for the discretization of (1.1) $A_h u_h = f_h$ amounts to the recursive application of the well-known V-cycle multigrid scheme (Algorithm 2.1). Here,

---
**Algorithm 2.1** Multigrid V-Cycle

| | | | |
|---|---|---|---|
| *Pre-smooth:* | $u_k$ | $\leftarrow$ | $S_k(u_k, f_k, A_k)$ |
| *Compute Residual:* | $r_k$ | $\leftarrow$ | $f_k - A_k u_k$ |
| *Restrict:* | $r_{k-1}$ | $\leftarrow$ | $R_k r_k$ |
| *Recurse:* | $e_{k-1}$ | $\leftarrow$ | $A_{k-1}^{-1} r_{k-1}$ |
| *Correct:* | $u_k$ | $\leftarrow$ | $u_k + P_k e_{k-1}$ |
| *Post-smooth:* | $u_k$ | $\leftarrow$ | $S_k(u_k, f_k, A_k)$ |

---

$S$ is the smoother and $k$ denotes the multigrid level. The solve at the coarsest level ($k = 0$), is done using a direct solver.

Multigrid methods can be classified into two categories, geometric and algebraic. The primary difference is that GMG methods rely on the underlying mesh connectivity for constructing coarser multigrid levels, whereas AMG methods are mesh-agnostic and work directly on the fine-grid matrix. The advantages of AMG are that it can be used as a black-box algorithm and does not require geometry or mesh information (only the assembled matrix). Its disadvantages are the communication-intense setup and the increased memory requirement compared to matrix-free geometric multigrid. Advantages of GMG are that it can be used in a matrix-free fashion and that it has low memory overhead. Moreover, operators can easily be modified at different levels, which can be necessary to accommodate certain boundary conditions. The disadvantages of GMG are that it requires a hierarchy of meshes and that it cannot be used in a black-box fashion.

Our GMG code is an extension of our previous work [39] to support high-order discretizations. For the AMG comparison, we use *ML* [21] from the *Trilinos* Project [24]. *ML* implements smoothed aggregation, a variant of AMG, and has shown excellent robustness and scalability. All experiments reported in this paper for multigrid (AMG and GMG) were performed using a single multigrid V-cycle as a preconditioner for the conjugate gradient (CG) method [44].

For the AMG and GMG scalability experiments, we report times (in seconds) for the following stages:
- *Setup*: setting up the multigrid hierarchy and computing the diagonal of the operator for the Jacobi smoother,
- *Eval*: applying the smoother and in the coarse grid solve, and
- *Comm*: performing restriction and prolongation.
- *Solve*: *Eval* + *Comm*.

In the rest of this section we provide additional details of our GMG implementation, in particular where it differs from [39].

**2.3.1. Meshing.** Our parallel geometric multigrid framework is based on hexahedral meshes derived from adaptive octrees [8, 37, 39]. As the multigrid hierarchy is independent of the discretization order, the construction of the grid hierarchy is identical to [39]. Multigrid requires the construction of a hierarchy of meshes, such that every element at level $k$, is either present at the coarser level $k-1$ or is replaced by an coarser(larger) element. The main steps in building the grid hierarchy are,

*Coarsening:.* Since we use a Morton-ordered linear representation for the octree, all eight sibling-octants, if present, will occur together. While partitioning the octants across processors, we ensure that all eight sibling-octants are owned by the same processor. Because the partitioning guarantees a single level of coarsening, the actual coarsening algorithm is embarrassingly parallel with $\mathcal{O}(N/p)$ parallel time complexity where $N$ is the total number of elements in the finer octree and $p$ is the number of processes. 2:1 balancing is enforced following the coarsening operation.

*Partition:.* Coarsening and the subsequent 2:1-balancing of the octree can result in a non-uniform distribution of coarse elements across the processes, leading to load imbalance. The Morton ordering enables us to equipartition the elements by performing a parallel scan on the number of elements on each process followed by point-to-point communication to redistribute the elements.

The reduction in the problem size at successive grids creates some problems from the perspective of load balancing. Since we partition each grid separately to ensure that all processes have an equal number of elements, we can end up with a very small number of elements per core (partitioning 2000 elements across 1000 cores) or even impossible cases (2000 elements across 10,000 cores). The surrogate mesh is generated (at every coarsening step) to facilitate the parallelization of intergrid transfers. While partitioning the surrogate mesh, we dynamically reduce the number of processes that are active at the coarser grid in order to ensure that communication does not dominate the computation. Additional details regarding our load balancing approach are discussed in [39].

*Meshing:.* By meshing we refer to the construction of the (numerical) data structures required for finite element computations from the (topological) octree data. In addition to the mesh extracted on the fine grid that is relevant for the simulation as a whole, we extract two meshes per multigrid level. First, we extract a *surrogate* mesh after coarsening and 2:1-balance of the fine mesh. Second, we extract the coarse mesh after repartitioning the surrogate mesh. This way we separate the processor-local nu-merical restriction (computation phase) from the parallel partition (communication phase). The fine and coarse meshes always contain all information for applying the elliptic operator $A_k$ in addition to the encoding of the partition. The surrogate mesh only contains the encoding of the partition and sizes of the elements. The use of the surrogate meshes for intergrid transfers is further detailed in [39].

*Recursion:.* To build the next multigrid level in the hierarchy, we repeat the previous steps using the coarse mesh as the new fine mesh. The recursion is stopped when either the required number of multigrid levels have been created or when no further coarsening of the mesh is possible.

**2.3.2. Discretization.** We use high-order discretizations based on Legendre-Gauss-Lobatto (LGL) nodal basis functions for polynomial orders $1 \leq q \leq 16$. For tensorized nodal basis functions on hexahedral meshes, the application of elemental matrices to vectors can be implemented efficiently by exploiting the tensor structure of the basis functions, as is common for spectral elements, e.g., [13]. This results in a computational complexity of $\mathcal{O}(q^4)$ for the element MatVec instead of $\mathcal{O}(q^6)$ if the element matrices were assembled. Globally, across $p$ processors, using tensor products the work complexity for a MatVec is $\mathcal{O}(Nq/p)$, requiring $\mathcal{O}(N)$ storage[4].

**2.3.3. Smoother.** We use damped Jacobi smoothing with $\omega = \frac{2}{3}$ for all runs reported in this paper. Although the performance of the Jacobi smoother deteriorates rapidly at $q > 4$ for variable coefficient problems, it performs well for the constant coefficient Poisson problem. For variable coefficient problems high-order Multigrid, Chebyshev-accelerated Jacobi smoother [7] provides good multigrid convergence up to $q = 16$. The cost of applying the Chebyshev-accelerated Jacobi smoother is similar to the Jacobi smoother, although it does require estimation of the maximum eigenvalues of the system matrix, which has to be done during the setup. For our experiments, we estimated the largest eigenvalue using 10 iterations of the Arnoldi algorithm. This increases the setup cost for multigrid.

**2.3.4. Restriction & Prolongation Operators.** Since the coarse-grid vector space is a subspace of the fine-grid vector space, any coarse-grid vector $v$ can be expanded independently in terms of the fine and coarse basis vectors,

$$(2.1) \qquad v = \sum_i v_{i,k-1} \phi_i^{k-1} = \sum_j v_{j,k} \phi_j^k,$$

where, $v_{i,k}$ and $v_{i,k-1}$ are the coefficients in the basis expansion for $v$ on the fine and coarse grids respectively.

The prolongation operator can be represented as a matrix-vector product (MatVec) with the input vector as the coarse grid nodal values and the output as the fine grid nodal values. The matrix entries are the coarse grid shape functions evaluated at the fine-grid vertices, $p_i$,

$$(2.2) \qquad P(i,j) = \phi_j^{k-1}(p_i).$$

The proof for (2.2) can be found in [37]. Similar to matrix-free applications of the system matrix, all operations required for the intergrid operations are done at the element level. As in [37], the restriction operator is the transpose of the prolongation operator.

---

[4]as opposed to $\mathcal{O}(Nq^3/p)$ work and $\mathcal{O}(Nq^3)$ storage

**2.3.5. Complexity.** For a $q$-order discretization with $N$ unknowns, the number of elements in the mesh is $\mathcal{O}(N/q^3)$. Given $p$ processes, the time complexity of building the multigrid hierarchy is $\mathcal{O}(N/(pq^3) + \log p)$. The $\mathcal{O}(N/(pq^3))$ corresponds to the coarsening and meshing operations, which are linear in the number of elements ($N/(pq^3)$). The second term accounts for the creation of the $\log p$ multigrid levels. The complexity of enforcing 2:1-balance is $\mathcal{O}(N/(pq^3) \log N/(pq^3))$ [9]. The cost of a `MatVec` is $\mathcal{O}(Nq/p)$. Assuming a uniform grid, we can estimate the communication costs as well: Each coarsening, balancing, and partition-correction call requires additional $\mathcal{O}(\log p)$ time to `MPI_Allgather` the local element count (one long integer). For partitioning and transferring, all elements of a process can potentially be communicated to $\mathcal{O}(1)$ processes. These transfers are implemented using non-blocking point-to-point communications, therefore the communication complexity per process is $\mathcal{O}(N/p)$. Thus the complexity for the GMG solve (without setup) is

$$T_{\text{GMG}} = \mathcal{O}\left(\frac{Nq}{p}\right) + \mathcal{O}\left(\log p\right).$$

**2.3.6. Caveats.** Although our GMG code is capable of handling complex geometries, the version used in this comparison is optimized for a unit cube domain with periodic boundary conditions. The setup phase, specifically the computation of the diagonal will be significantly more expensive for other domains. In addition, the smoother will also be more expensive due the need to compute the Jacobian of the mapping (this can be traded for pre-computation and additional storage) from the reference element to each element.

**3. Experimental setup.** In this section, we give details on the experimental setup we used to test the methods.

*Hardware.* The hardware employed for the runtime experiments carried out is the Stampede system at TACC and Titan at ORNL. Stampede entered production in January 2013 and is a high-performance Linux cluster consisting of 6,400 compute nodes, each with dual, eight-core processors for a total of 102,400

CPU-cores. The dual-CPUs in each host are Intel Xeon E5-2680 2.7GHz (Sandy Bridge) with 2GB/core of memory and a three-level cache. The nodes also feature the new Intel Xeon Phi co-processors. Stampede has a 56GB/s FDR Mellanox InfiniBand network connected in a fat tree configuration which carries all high-speed traffic (including both MPI and parallel file-system data). Titan is a Cray XK7 with a total of 18,688 nodes consisting of a single 16-core AMD Opteron 6274 processor, for a total of 299,008 cores. Each node has 32GB of memory. It is also equipped with a Gemini interconnect and 600 terabytes of memory across all nodes.

*Experiments.* We study cost per unknown and cost per accuracy for two sets of experiments:

- In the first set of tests we consider the ***cost per unknown for a fixed algebraic accuracy*** for all five different methods, both in terms of computation and communication. This is ***not*** equivalent to which method is faster. It just reveals the constants in the complexity estimates, and the effect of communication and overheads. For example, high-order methods have better locality but higher work per unknown, so it is not entirely clear what the effect is in the work per unknown. Also, we only test regular grids, but *GMG and FMM have overheads since they use pointer-based data-structure to support octree discretizations.*

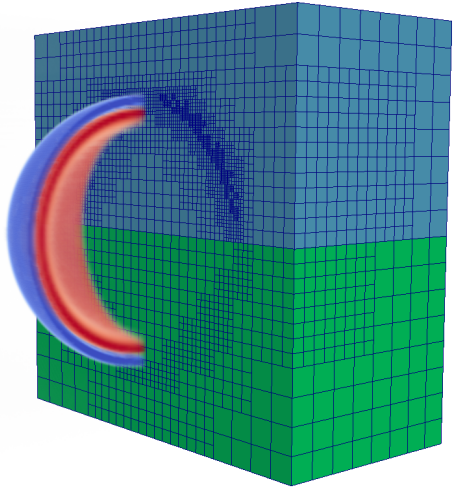

**Fig. 3:** *Adaptive mesh structure for the forcing term, f, in test case 2. The green cube (bottom quadrant) is the adaptive mesh using 6th order elements and the blue cube (top quadrant) is that of using 14th order elements. The number of unknowns is significantly reduced.*

In Figures 4 and 5 and Tables 3, 4, 5, and 6, we consider weak and strong scalability of the solvers with regard only the application of the "inverse". For all cases we use an arbitrary right-hand side.

For the FFT, it means an 3D FFT , diagonal scaling to apply the inverse of the Laplacian in the spectral domain, and then another 3DFFT. For the FMM, is a an application of the scheme to the right hand side. For the AMG and GMG runs, we use CG preconditioned with multigrid and we drive the relative residual down to 1E-13. We use the V-cycle (using a full cycle makes little difference). In these runs, on average GMG-1 takes 6 iterations, GMG-4 takes 8 iterations, GMG-8 takes 20 iterations, and GMG-16 takes 28 iterations. These are CG iterations in which the preconditioner is the corresponding multigrid method. For the hpGMG runs we just used the default settings and just controlled the tolerance.

To give an estimate of the relative constants and connect to the complexity estimates we compile statistics for all runs on Stampede and we report the results in Table 11.

- In the second set of experiments, **we examine the time to solution for the different methods**, irrespective of order, or grid. We just set the number of processors and the target accuracy. We report timings for solution error that is 1E-7 or less. We report timings for two different exact solutions $u$, one that requires uniform discretization (Table 7) and one that requires a highly refined discretization (Table 8). For the latter, we also consider lower accuracy (Table 9) and we also compare with the AMG solver (Table 10).

  The first exact $u$ is an oscillatory field given by $u(x_1, x_2, x_3) = \sin(2\pi k x_1)\sin(2\pi k x_2)\sin(2\pi k x_3)$. This field requires uniform discretization. FFT can resolve $u$ to machine precision using $2k$ points per dimension. First-order methods (for example the 7-point Laplacian or linear finite elements) require roughly $4k$ points per digit of accuracy. Roughly speaking, for two digits of accuracy in 3D, first-order methods require $64\times$ more unknowns than FFT, and for seven digits of accuracy require $2048^3 k$ more unknowns than FFT. Higher-order methods, reduce the number of unknowns significantly but are still suboptimal compared to FFT for this particular $u$.

The second exact $u$ is $\exp(-(r/R)^\alpha)$, which has a sharp internal layer at the surface of a sphere of radius $R$, where $r$ is the distance from the point we evaluate $u$ to the center of the sphere. As $\alpha$ increases $u$ develops very sharp gradients around $r = R$ and is almost constant everywhere else. Roughly speaking, every time we double $\alpha$ we increase the spectrum of $u$ by a factor of two and therefore in 3D the number of unknowns for FFT increases by a factor of eight. On the other hand, an adaptive method only refines around the area of the internal layer and can resolve the solution for large $\alpha$. Figure 3 shows an example of an adaptive mesh used by FMM. As a result, it is expected that either FMM or GMG/AMG will be the optimal method for this test case.

For both time-to-solution test cases, the source function is set to the closed form Laplacian of the exact solution, that is $f = -\Delta u$. This $f$ is then used to compute the discretized solution $u_N$. We report $\|u - u_N\|_{\ell_\infty}$, where we abuse the notation and use $u$ for the exact solution evaluated at the discretization points.

*Software parameters.* For all our runs FFT means we used AccFFT. For some of the larger Titan runs we have used P3DFFT. For the FMM we have two settings, high accuracy ($q = 14$, $m = 10^2$) and low accuracy $q = 6$ and $m = 4^2$. For FMM the adaptivity criterion is based on the decay of the tail of the Chebyshev expansion of $f$ at every leaf node. For GMG and AMG the adaptivity is based on the gradient of $f$. For the multigrid case the smoothing is done using a pointwise Jacobi scheme with two pre-smoothing steps and one post-smoothing step. The relative algebraic residual is driven (using CG) to 1E-13 in all runs. The coarse solves were done using sparse parallel LU factorization (SuperLU [28]).[5]

**4. Results.** We are evaluating our geometric multigrid (**GMG**) with $q = 1, 4, 8, 10, 14$, our fast multipole method (**FMM**) with $q = 6, 14$, the **hpGMG** ($q = 1$ and $q = 2$), and finally the (**FFT**) using the AccFFT library. As a reference we also perform some runs with Trilinos's ML algebraic multigrid (**AMG**) library. When using AMG, we are using matrices that are generated by our library using first and fourth order polynomials. AMG is used as a reference point to show the difference between using a state of the art general-purpose solver over a more problem-specific solver. All solvers use only CPU cores; the Phi coprocessors were not used.

We use the following notation to denote the different variants of the methods we are testing: "$X$-$q$" indicates method "$X$" ran with $q$-th order polynomials. For example "GMG-4" indicates geometric multigrid with fourth order polynomials; "AMG-1", indicates algebraic multigrid with linear polynomials.

Our first test is a weak and strong scaling test on Stampede and Titan to assess the cost per unknown as a function of the method and the problem size without regard to the accuracy. The strong scaling results are reported in Tables 3 (Stampede) and 5 (Titan) and the weak scaling results are reported in Tables 4 (Stampede) and 6 (Titan). All these results where done using uniform discretization so that we could compare with FFT. FMM and GMG have overhead costs since they are not optimized specifically for regular grids.

---

[5] The main parameters for ML are aggregation threshold 0.2, aggregation type : uncoupled, smoother Gauss Seidel (three iterations), and the coarse grid solver is KLU.

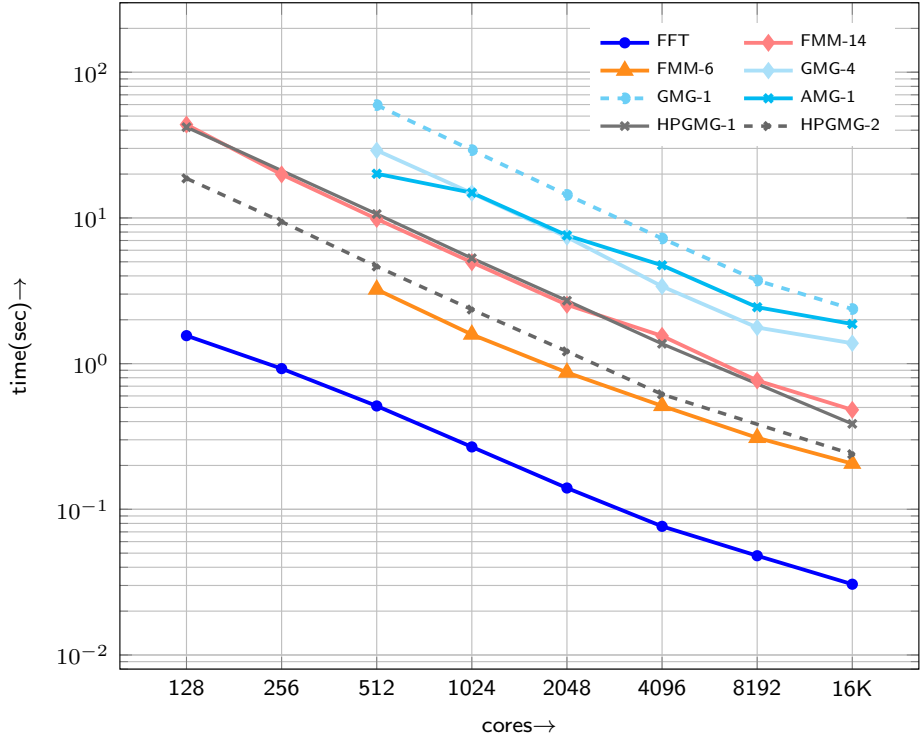

**Fig. 4:** *Here we report strong scaling results on Stampede that capture the cost per unknown for the different methods. These are only the "Solve" times compared to Table 3. For the FFT and FMM this involves just a single evaluation. For the GMG and AMG methods it involves several multigrid sweeps to bring the relative algebraic residual to zero (up to double precision). All times are in seconds for a problem with $N = 1024^3$. The parallel efficiency for the different methods is: FFT (52%), FMM-14 (64%), AMG-1 (33%), GMG-4 (67%), HPGMG-1 (96%), and HPGMG-2(94%). In absolute terms FFT has always the lowest constant per unknown. We can observe orders of magnitude differences between the different variants. Since FMM involves just an evaluation lower order has a smaller constant than higher. Notice also that for our implementations, FMM with high order is significantly faster than GMG with low order. Although, we don't report the results in detail here, we also tested AMG-4 and it was 2× slower than AMG-1. Note that HPGMG is specifically optimized for uniform discretization, and thus is faster compared to our GMG library which can also handle adaptive high-order multigrid.*

The second set of experiments examine the time to solution, which are of course the ones that fully characterize the performance of these solvers. The results for the oscillatory $u$ are reported in Table 7 and the results for the non-uniform $u$ are reported in Table 8.

We also report wall-clock times of "*Setup*", "*Comm*" and "*Solve*". The "*Setup*" involves costs that incur while building data-structures whenever a change in the number of unknowns or their spatial distribution takes place. For the GMG, AMG, building the grid hierarchy, 2:1 balancing it, computing coefficients for the Chebyshev smoother, and setting-up the sparse direct solver for the coarse solve. For the FMM, it involves loading precomputed near-interaction translation operators from the disk to main memory, computing far-field interaction operators, allocating memory buffers, building the octree and its 2:1 balance refinement. For right-hand sides that require the same number of unknowns but different trees due to different distribution of spatial features, the setup has to be called again. For the FFT the setup phase involves parameter tuning to optimize CPU performance. "*Solve*" is the time required to evaluate the solution $u$ whenever a new right-hand side $f$ is specified (without changing the resolution); and "*Comm*" is the distributed memory communication costs during the "*Solve*" phase. For AMG library we do not report communication costs, as we could not find an easy way to measure it (as of version ML 5.0). As we mentioned in the introduction more detailed analysis of the three algorithms (along with discussion on single-core performance) can be found in [30, 33, 39].

| $p$ | phase | AMG-1 | GMG-1 | GMG-4 | HPGMG-1 | HPGMG-2 | FMM-14 | FFT |
|---|---|---|---|---|---|---|---|---|
| 1024 | Setup | 5.4 | 1.5 | 0.3 | - | - | 5.2 | 5.9 |
| | Comm | - | 2.8 | 2.9 | - | - | 0.2 | 0.1 |
| | Solve | 15 | 26.2 | 11.8 | 5.3 | 2.3 | 5.0 | 0.3 |
| 4096 | Setup | 4.1 | 0.7 | 0.3 | - | - | 3.9 | 2.5 |
| | Comm | - | 0.9 | 0.9 | - | - | 0.1 | 0.05 |
| | Solve | 4.7 | 6.3 | 2.5 | 1.4 | 0.6 | 1.5 | 0.08 |
| 16384 | Setup | 15 | 0.5 | 0.4 | - | - | 5.2 | 1.4 |
| | Comm | - | 0.9 | 0.8 | - | - | 0.1 | 0.02 |
| | Solve | 1.9 | 1.6 | 0.6 | 0.4 | 0.2 | 0.5 | 0.03 |

**Table 3:** *Here we report strong scaling results on Stampede. We report representative breakdowns of the timings (in seconds) for uniform grids with $N = 1024^3$. GMG has insignificant setup costs. AMG has a higher cost as it needs to build the multigrid hierarchy algebraically using graph coarsening. Finally, FMM has setup costs that are in the order of the solve time; these costs are related to tree construction, loading (from the disk) precomputed tables, 2:1 balancing, and constructing the interaction lists. The precomputed tables are generic and valid for both uniform and non-uniform grids. The setup phase for FMM has fixed size overheads associated with loading and computing interaction operators and therefore, the setup time does not scale as we increase the number of processor cores. Note that HPGMG is specifically optimized for uniform discretization, and thus is faster compared to our GMG library which can also handle adaptive high-order multigrid.*

| $p$ | phase | AMG-1 | GMG-1 | GMG-4 | HPGMG-1 | HPGMG-2 | FMM-14 | FFT |
|---|---|---|---|---|---|---|---|---|
| 128 | Setup | 4.9 | 1.7 | 0.1 | - | - | 4.1 | 5.0 |
| | Comm | – | 2.0 | 2.1 | - | - | 0.2 | 0.1 |
| | Solve | 14.9 | 32.9 | 15.8 | 5.4 | 2.4 | 4.9 | 0.2 |
| 1024 | Setup | 5.4 | 1.5 | 0.3 | - | - | 4.5 | 5.9 |
| | Comm | – | 2.8 | 3.0 | - | - | 0.2 | 0.16 |
| | Solve | 14.9 | 32.2 | 16.9 | 5.3 | 2.3 | 5.1 | 0.26 |
| 8192 | Setup | 6.5 | 1.7 | 0.5 | - | - | 5.6 | 8.0 |
| | Comm | – | 4.3 | 4.6 | - | - | 0.3 | 0.22 |
| | Solve | 15.2 | 35.7 | 17.9 | 7.6 | 3.4 | 5.3 | 0.34 |

**Table 4:** *Here we report weak scaling results on the Stampede system. In this test, we keep the number of unknowns per process fixed as we increase $p$. The problem sizes that used are $N = 512^3$, $1024^3$, $2048^3$ for $p = 128$, $1024$, $2048$ respectively. For the AMG and HPGMG runs, we could not directly measure the communication costs.*

*Strong scaling analysis on Stampede, uniform grid.* To demonstrate the overhead and performance characteristics of the methods, we consider strong scaling for a problem with one billion unknowns ($1024^3$). We report "Solve" time in Figure 4; and the "Setup" and "Comm" times in Table 3. We report "Solve" time separately since in many cases setup time is amortized across nonlinear iterations or time-stepping. As one can see, FFT is clearly the fastest method, although FMM is not far behind. Here AMG-1, despite having less computation per unknown is somewhat slower due to the cost per V-cycle and the number of sweeps. All the methods scale quite well. Also notice that GMG-4 is faster than GMG-1 due to better compute intensity. Same is true for HPGMG-2 compared to HPGMG-1. Note that HPGMG is specifically optimized for uniform discretization, and thus is faster compared to our GMG library which can also handle adaptive high-order multigrid. The FFT setup time only depend on the problem size and number of processors and thus is done just once during the course of a calculation. In most applications it is unusual to change resolution dynamically when using FFTs.

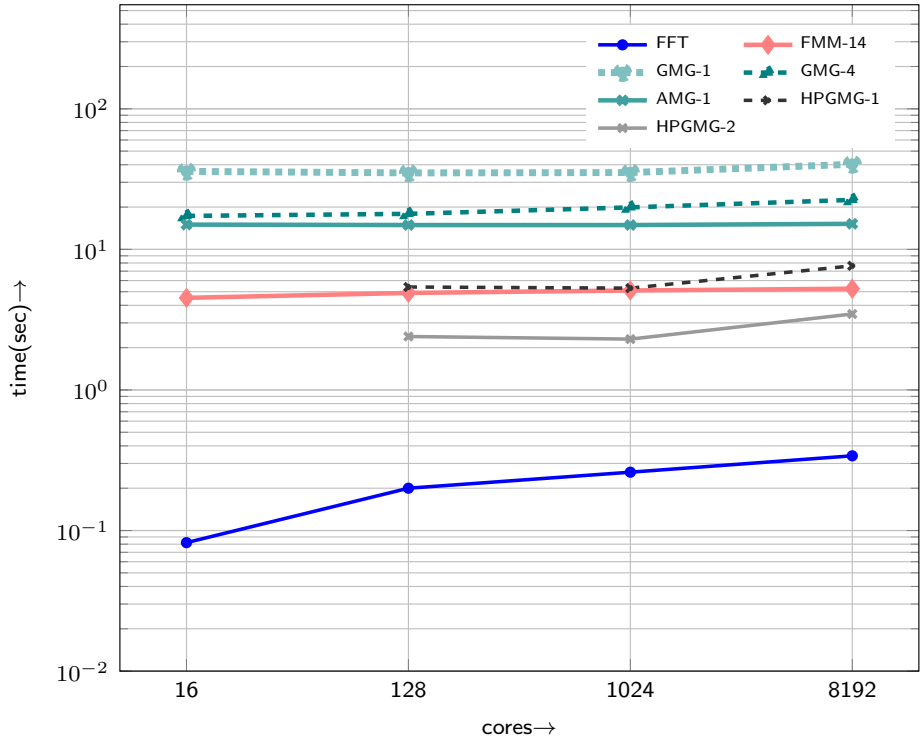

**Fig. 5:** *Weak scaling on Stampede: We report "Solve" time in seconds of different methods for a problem with a grain size of one million unknowns per core; the largest problem on Stampede has 8.5 billion unknowns on 8,192 cores. All the methods scale quite well. FFT is the fastest and the difference is over 30× over GMG-1.*

| $p$ | phase | FMM-14 | P3DFFT |
|---|---|---|---|
| 57,344 | Setup | 105 | — |
|  | Comm | 14.4 | — |
|  | Solve | 230 | 44.2 |
| 114,688 | Setup | 44.4 | — |
|  | Comm | 6.3 | — |
|  | Solve | 105 | 19.1 |
| 229,376 | Setup | 30 | — |
|  | Comm | 4.2 | — |
|  | Solve | 58 | 11.3 |

**Table 5:** *Strong scaling on Titan (for FMM only). We report wall-clock time in seconds for FMM with 14th order polynomials for a problem with 600 billion unknowns. The efficiency is 99% for the FMM "Solve" phase. The communication costs for the FMM are order-of-magnitude smaller than the "Solve" and "Setup" phases.*

*Weak scaling on Stampede.* To illustrate the constant factors in the complexity estimates for the different methods, we consider a weak scaling test in which we keep the number of unknowns per core fixed to roughly one million and we increase $N$ and $p$, keeping their ratio fixed Figure 5. The observed end-to-end efficiencies for the different methods are as follows: FFT(59%), FMM-14(95%), GMG-1(92%), GMG-4(88%), AMG-1 (98%), HPGMG-1 (69%), and HPGMG-2 (70%). FFT does not scale as well as other methods but it is still quite faster. It is worth noting that although GMG-4 does not scale as well as GMG-1, but it is about two times faster. This is due to the higher order scheme used in GMG-4.

| $N$ | $p$ | phase | FMM-14 | AccFFT | P3DFFT |
|---|---|---|---|---|---|
| $1,024^3$ | 448 | Setup | 16.0 | 12.0 | — |
| | | Comm | 2.5 | 1.3 | — |
| | | Solve | 51.8 | 2.2 | 4.5 |
| $2,048^3$ | 3,584 | Setup | 16.7 | 15.0 | — |
| | | Comm | 2.7 | 1.9 | — |
| | | Solve | 52.5 | 2.6 | 6.0 |
| $4,096^3$ | 28,672 | Setup | 17.6 | — | — |
| | | Comm | 3.6 | — | — |
| | | Solve | 53.1 | — | 8.0 |
| $8,192^3$ | 229,376 | Setup | 19.0 | — | — |
| | | Comm | 4.2 | — | — |
| | | Solve | 57.8 | — | 11.3 |

**Table 6:** *Weak scaling results on Titan. We report wall-clock time in seconds for the FFT and FMM with 14th order polynomials. The largest problem size corresponds to half a trillion unknowns. FFT is significantly faster than FMM.*

*Strong and weak scaling on Titan.* On Titan we have performed weak scaling with FFT, and FMM-14. The low-order versions for these methods have already scaled to 260K cores on Jaguar [39].)

The strong scaling results are summarized in Table 5 for a problem with 600 billion unknowns for the FMM.

The weak scaling results for FFT and FMM are shown in Table 6. We observe that the FMM is almost 25× slower than the FFT. This is because FFT has much less work per unknown. Also it is "more" memory bound and less sensitive to lower per-node peak performance (which is the case for Titan CPUs compared to Stampede CPUs) so FMM runs slower on Titan whereas FFT is not affected as much.

*Time to solution weak scaling.* So as we saw in the previous paragraphs, one can test these solvers using the same number of unknowns across solvers, where significant differences were observed. However, what matters most is time to solution. Once this is factored in, the differences between the solvers become more pronounced.

We start by looking at test case 1 (the oscillatory synthetic problem), a problem that requires uniform refinement. We report the "Solve" time to reach single precision accuracy in Table 7. From this first experiment, it is clear that FFT is extremely efficient for this problem. For example, for an 8,192-core run, if we retain the same grid and we change $k$ from 128 to 256, the errors for FMM and GMG drop to three digits of accuracy whereas FFT can still resolve it to machine precision (as a reminder $k$ is the frequency of the $f$). This is of course expected given the spectral accuracy of the FFT and the very low work per unknown required by the algorithm. The purpose of the experiment is to highlight the performance gap even if we use very high order FMM or Galerkin methods.

But what about solutions with localized features as in test case 2? We report the "Solve" time $u = \exp(-(r/R)^\alpha)$ for Table 8 as a function of $\alpha$. The larger the $\alpha$ the more localized features $u$ has. The performance of the different solves is now quite different. FFT works well for relatively small values of $\alpha$ but as we make the solution sharper (larger $\alpha$) it cannot resolve the length-scales of the solution and the uniform grid size shoots up; for the largest problem size FFT requires 100× more unknowns than those required by FMM and GMG. Furthermore, FMM outperforms multigrid but the difference is not as dramatic. Understanding this requires a lengthier analysis and is sensitive to the quality of per core performance. In passing, let us mention that the adaptive multigrid is quite efficiently implemented. Both codes used dense matrix-matrix multiplication kernels (that is they used BLAS DGEMM routine) but their intensity is quite different. GMG uses a high order tensor basis. The main computational kernel for GMG is the element traversal to apply the local stiffness matrix to a vector. This kernel is decomposed to 1D slices to reduce the complexity. The decomposition involves more integer operators and multiple DGEMM invocations. For near interactions FMM truncates the tensor basis and uses a single DGEMM of larger size than GMG for the near interaction. The far interaction is much more involved and briefly summarized in §2.2.

|   |   | **GMG-14** | | | **FMM-14** | | | **FFT** | |
|---|---|---|---|---|---|---|---|---|---|
| $p$ | k | T | $\ell_\infty$ | L | T | $\ell_\infty$ | L | T | L |
| 16 | 16 | 10 | 1E-8 | 4 | 4.5 | 3E-7 | 5 | <0.01 | 6 |
| 128 | 32 | 9.8 | 4E-8 | 5 | 4.9 | 3E-7 | 6 | <0.01 | 7 |
| 1,024 | 64 | 10 | 1E-7 | 6 | 4.9 | 4E-7 | 7 | <0.01 | 8 |
| 8,192 | 128 | 10 | 3E-7 | 7 | 5.1 | 5E-7 | 8 | <0.01 | 9 |

**Table 7:** *Time-to-solution for test case 1 (oscillatory field) performed on Stampede: Here we report the number of cores p, the effective resolution (wavenumber) k, the relative discrete infinity norm $\ell_\infty$ of the error in u, "Solve" time "T" and the uniform refinement level "L", required to achieve single-precision (seven digits) accuracy or better. For the FFT we do not report errors since it resolves the solution to machine accuracy. FFT dramatically outperforms FMM and GMG, since it can resolve the problem with considerably fewer unknowns. FFT requires $8^L$ unknowns to resolve f. FMM requires $8^L q^3/6$ (and GMG $8^L q^3$) or roughly 9 billion unknowns to resolve a u with k = 128 to single precision; FFT can resolve k = 256 to machine precision with nearly 70× fewer points. Also note that the target precision doesn't change the conclusions, it just determines the overall problem size.*

|   |   | **GMG-5** | | **FMM-14** | | **FFT** | |
|---|---|---|---|---|---|---|---|
| $p$ | $\alpha$ | T | N | T | N | T | N |
| 32 | 10 | 4.4 | 4.4E+6 | 0.5 | 2.9E+6 | 0.1 | 1.7E+7 |
| 512 | 40 | 5.0 | 1.8E+7 | 0.7 | 6.4E+7 | 0.5 | 1.1E+9 |
| 8,192 | 160 | 6.0 | 2.5E+8 | 0.8 | 1.1E+9 | 2.6 | 6.9E+10 |
| 32,768 | 320 | 6.6 | 9.8E+9 | 1.8 | 4.3E+9 | 4.9 | 5.5E+11 |

**Table 8:** *Time-to-solution for test case 2 performed on Stampede: Here $\alpha$ is a measure of the difficulty of resolving u. The higher the $\alpha$ value the sharper the derivatives. We report "Solve" time "T" and the required number of unknowns N to achieve single precision accuracy or higher. In all of these runs and for all three methods the relative error $\ell_\infty$ is less than $4E-7$. We used GMG-5, because in our test it was the fastest variant of GMG (see Tab. 10).*

Due to the significance difference between the GMG and FMM kernels, the result is that their floating point performance per core differs. For our implementations FMM outperforms GMG. We also report AMG/GMG timings for different orders in Table 10. As expected AMG is significantly slower than GMG. We can see that the time to solution decreases by increasing the order up to $q = 5$, and then increases afterwards. Finally, we would like to remark that both GMG and FMM over-refine to maintain 2:1 balancing.

It would be tempting to use the expressions for $T_{\text{FMM}}, T_{\text{FFT}}, T_{\text{GMG}}$ of §2 with N fixed to try to quantify the effect of the order q on the communication and computation costs and how it compares between the methods. However N *does* depend on q and the relation is non-trivial. That's why our empirical results paint a more accurate

|       | L  | N      | $\ell_\infty$ | T    |
|-------|----|--------|---------------|------|
| GMG-1 | 12 | 9.5E+6 | 9.4E-3        | 9.40 |
| GMG-2 | 9  | 1.0E+6 | 5.3E-3        | 1.60 |
| GMG-5 | 6  | 4.1E+5 | 2.2E-3        | 3.78 |
| FMM-7 | 8  | 2.7E+6 | 3.5E-3        | 0.13 |
| FFT   | –  | 5.6E+7 | 8.1E-3        | 0.07 |

**Table 9:** *Time-to-solution for test case 2, using 32 cores on Stampede. We set $\alpha = 40$. We report "Solve" time "T" and the required number of unknowns N and tree refinement level "L" to achieve about two digits of accuracy or higher. For the multigrid-preconditioned runs, CG took 2,3, and 5 iterations for the V-cycle GMG-1, GMG-2, and GMG-5 schemes.*

| p  | q | N     | T (AMG) | T (GMG) |
|----|---|-------|---------|---------|
| 32 | 4 | 8.9E6 | 29.3    | 8.6     |
| 32 | 5 | 4.4E6 | 27.7    | 4.4     |
| 32 | 6 | 2.1E6 | 22.8    | 4.9     |
| 32 | 8 | 6.0E5 | 51.4    | 6.1     |

**Table 10:** *Corresponding results in Tab. 8 with $\alpha = 10$ using AMG/GMG with different orders. As the order is increased up to 6th order, the cost per unknown reduces. However, further increase in the order increases this cost and thus the time to solution increases.*

picture.

*Performance of GMG.* Our comparisons against the finite element version of hpGMG Table 3 indicate that hpGMG performed $4-5\times$ faster than GMG-1. While the performance of hpGMG is impressive and reflective of its quality as a benchmark, our performance is equally strong given that our code supports adaptivity and arbitrarily high orders; hpGMG only supports linear and quadratic elements. We briefly summarize the main reasons for the $\sim 5\times$ speedup observed with hpGMG.

- Given that our code supports adaptive meshes, we need to store additional mesh-lookup information. This reduces the effective flop-to-mop ratio of our code.
- Adaptivity also forces us to have to deal with nodes that do not represent independent degrees of freedom. This increases the cost of the elemental matrix-vector multiplications. The details of how we treat such nodes is beyond the scope of this paper. Please refer to [37, 39] for two different approaches to handling these cases.
- hpGMG uses a cartesian grid topology, leading to efficient communication. Our GMG code is capable of handling complex geometries in an adaptive fashion and we use a Morton-ordering for partitioning work across processors. This makes it difficult for us to provide any topological information to enable efficient communication.

We also compared the performance of hpGMG for the oscillatory synthetic problem reported in Table 7. For the case of $p, k = 16$, using hpGMG with $Q1$ elements, we could accomodate a mesh of size $512^3$ requiring 42.3 seconds to reach a relative $\ell_\infty$ norm of the error in u of 1.1E-2. Using $Q2$ elements, we could accomodate $256^3$ elements, to reduce the error to 2.8E-3 in 19 seconds. Assuming perfect convergence rate (error = error/8 for each refinement), we will need $4096^3$ elements to reduce the error to 7E-7. Assuming linear scaling, we will need $77K$ seconds, compared to us needing 10 seconds. This example illustrates that performance uniform grid discretizations deteriorates significantly for problems that have layers and localized features.

| Method | Communication | | Computation | |
|---|---|---|---|---|
| | $T_{\text{COMM}}$ | $\sigma_{\text{COMM}}$ | $T_{\text{COMP}}$ | $\sigma_{\text{COMP}}$ |
| FFT | $1.7\text{E-}7\,\frac{N}{p}$ | 7.9E-8 | $4.2\text{E-}9\frac{N\log N}{p}$ | 1.2E-9 |
| FMM-14 | $5.3\text{E-}6\left(\frac{N}{p}\right)^{2/3}\log p$ | 2.0E-6 | $4.5\text{E-}6\frac{N}{p}$ | 5.8E-7 |
| GMG-1 | $3.2\text{E-}4\left(\frac{N}{p}\right)^{2/3}$ | 1.3E-4 | $2.8\text{E-}5\frac{N}{p}$ | 3.9E-6 |
| GMG-4 | $3.2\text{E-}4\left(\frac{N}{p}\right)^{2/3}$ | 1.1E-4 | $1.3\text{E-}5\frac{N}{p}$ | 3.2E-6 |

**Table 11:** Here we report the constants of the complexity estimates of each method, for the communication time $T_{\text{COMM}}$ and the computation time $T_{\text{COMP}}$. The constants were computed using the data of Table 3 and 4. We also report the standard deviation ($\sigma_{\text{COMM}}$, $\sigma_{\text{COMP}}$) for the constants.

**5. Conclusions.** We outlined three solvers developed in our group for Poisson problems with high-order discretizations and we compared them for the Poisson problem with constant coefficients on the unit cube with periodic boundary conditions. To measure the efficiency of our schemes with the current state of the art we also compared our GMGM to implementations developed by other groups, in particular the finite-element HPGMG and AMG. There is no (parallel) software available that compares to our FMM.

- If we fix the number of unknowns, and we consider just the cost per unknown, FFT is the fastest scheme because it has really small constants in its complexity estimates. FMM performs comparably with FFT for the same number of unknowns (on Stampede). Which method is preferable depends on the type of resolution scheme (uniform vs non-uniform). For this reason, one should note that *it is not possible to compare these methods using exactly the same N*. Often in the literature, in particular the HPC literature, such comparisons appear. Without accuracy for a specific problem such information can be misleading.
- FFT outperforms the other methods by very large margins for uniform grids.[6] As expected the communication cost of FFT does not scale as favorably. However, FFT requires far fewer unknowns for functions that are within its resolution, and as we mentioned it also has the smallest constants. So even if the communication scaling is suboptimal, overall the method is much faster. We report results up to 230K cores that confirm this.
- For fields that have sharp local features regular-grid discretizations lose their advantage. Adaptively refined algorithms like GMG and FMM can resolve such fields more effectively since they require orders of magnitude fewer unknowns than FFT.
- GMG, hpGMG, or AMG based on linear elements are significantly slower compared to high-order methods when the accuracy requirements are high. The excessive number of unknowns can result in orders of magnitude in performance penalties. We expect this to hold for problems for which the solution has localized discontinuities.
- For non-uniform grids, both FMM and GMG methods scale quite well since they are matrix free and both use octree-based hierarchical space decompo-sition methods. The structure of the calculation however is different, with the FMM offering significant opportunities for task parallelism and compute intensity. On the other hand, the setup of FMM is more expensive. Our GMG is much more general than the FMM since it supports more general geometries and boundary conditions. All these incur overheads.
- A state of the art algebraic solver is orders of magnitude slower than the fastest variants of FFT, FMM and GMG. Of course, AMG is much more general as it can deal with arbitrary geometries, variable and anisotropic coefficients, and more general boundary conditions. However, many important applications require fast Poisson solvers on a cube with boundary conditions that FFT, FMM, or GMG can handle. Therefore, using a method designed for more general situations can incur quite high performance penalties.
- All three methods can be further improved: by integrating tightly with accelerators, further improving single-node performance, and trying to improve the per leaf high-order calculations (for FMM and GMG).
- No method rules the others. All methods have regimes in which they are the fastest. A robust black-box solver should be a hybrid method for cases when the true solution is a superposition of a highly oscillatory field and a highly localized field. This is specially true for problems with more complex geometries and boundary conditions.

---
[6]Note that the hpGMG benchmark includes a finite-volume stencil code for constant coefficients that it is extremely fast and possibly compares well with FFT. Also it is more general since it supports different types of boundary conditions. So our conclusion regarding FFT pertains strictly to the packages we used and it is not general.

**Acknowledgment.** The authors would like to thank Omar Ghattas and Georg Stadler for useful discussions on this topic. Also, we are grateful to ORNL and TACC for help with the experiments on the Titan and Stampede systems. Support for this work was provided by: the U.S. Air Force Office of Scientific Research (AFOSR) Computational Mathematics program under award number FA9550-09-1-0608; the U.S. Department of Energy Office of Science (DOE-SC), Advanced Scientific Computing Research (ASCR), DE-FG02-09ER25914, and the U.S. National Science Foundation (NSF) Cyber-enabled Discovery and Innovation (CDI) program under awards CMS-1028889 and OPP-0941678, and the PetaApps program under award OCI-0749334. Computing time on the Cray XK7 Titan was provided by the Oak Ridge Leadership Computing Facility at Oak Ridge National Laboratory, which is supported by the Office of Science of the Department of Energy under Contract DE-AC05-00OR22725. Computing time on the Texas Advanced Computing Centers Stampede system was provided by an allocation from TACC and the NSF.